\documentstyle[12pt,graphicx,amssymb]{article}
\begin{document}
\date{}
\renewcommand{\thefootnote}{}

\def\thebibliography#1{\noindent{\normalsize\bf References}
 \list{{\bf
\arabic{enumi}}.}{\settowidth\labelwidth{[#1]}\leftmargin\labelwidth
 \advance\leftmargin\labelsep
 \usecounter{enumi}}
 \def\newblock{\hskip .11em plus .33em minus .07em}
 \sloppy\clubpenalty4000\widowpenalty4000
 \sfcode`\.=1000\relax}

\def\fin{\hspace{5mm}$\square$}
\def\fig#1{\begin{center}{\fbox{\ Figure\ #1 \ }}\\[4mm]\end{center}}

\def\figure#1#2#3{
\begin{center}
\includegraphics[trim=0mm 0mm 0mm 0mm, width=.#1\linewidth]
{vassi#2.eps}\\
Figure #3
\end{center}}

\title{\vspace*{-1cm}
{\large LOCAL MOVES ON SPATIAL GRAPHS \\
AND FINITE TYPE INVARIANTS}}

\author{{\normalsize KOUKI TANIYAMA}\\
{\small Department of Mathematics, School of Education,
Waseda Universitry}\\[-1mm]
{\small Nishi-Waseda 1-6-1, Shinjuku-ku, Tokyo 169-8050, Japan}\\[-1mm]
{\small e-mail: taniyama@mn.waseda.ac.jp}\\[3mm]
{\normalsize AKIRA YASUHARA }\\
{\small Department of Mathematics, Tokyo Gakugei University}\\[-1mm]
{\small Nukuikita 4-1-1, Koganei, Tokyo 184-8501, Japan}\\
{\small {\em Current address}, October 1, 1999 to September 30, 2001:}\\[-1mm]
{\small Department of Mathematics, The George Washington University}\\[-1mm]
{\small Washington, DC 20052, USA}\\[-1mm]
{\small e-mail: yasuhara@u-gakugei.ac.jp}
}

\maketitle

\vspace*{-5mm}
\baselineskip=16pt
{\small
\begin{quote}
\begin{center}A{\sc bstract}\end{center}
\hspace*{1em}
We define $A_k$-moves for embeddings of a
finite graph into the 3-sphere for each natural number $k$.
Let $A_k$-equivalence denote an equivalence relation
generated by $A_k$-moves and ambient isotopy. $A_k$-equivalence implies
$A_{k-1}$-equivalence. Let ${\cal F}$ be an $A_{k-1}$-equivalence class
of the embeddings of a finite graph into the 3-sphere. Let ${\cal G}$ be
the quotient set of ${\cal F}$ under $A_k$-equivalence.
We show that the set ${\cal G}$ forms an abelian group under a certain
geometric operation. We define finite type invariants on ${\cal F}$ of
order $(n;k)$.  And we show that if any finite type invariant of order
$(1;k)$ takes the same value on two elements of ${\cal F}$, then they
are $A_k$-equivalent.
$A_k$-move is a generalization of $C_k$-move defined by K. Habiro.
Habiro showed that two oriented knots are the same up to $C_k$-move and
ambient isotopy if and only if any Vassiliev invariant of order $\leq k-1$
takes the same value on them. The \lq if' part does not hold for
two-component links. Our result gives a sufficient condition for
spatial graphs to be $C_k$-equivalent.
\end{quote}
}

\footnote{{\em Short Running Title}: Local Moves on Spatial
Graphs and Finite Type Invariants}

\newpage

\baselineskip=18pt 

\noindent
{\bf Introduction}

\bigskip
K. Habiro defined a local move, {\it $C_k$-move},
for each natural number $k$ \cite{Habiro2}.
It is known that if two embeddings $f$ and $g$ of a graph into
the three sphere are the same up to $C_k$-move and ambient isotopy,
then $g$ can be deformed into a band
sum of $f$ with certain $(k+1)$-component links and
that changing position of a band and an arc, which is
called {\em a band trivialization} of $C_k$-move, is realized
by $C_{k+1}$-moves and ambient isotopy \cite{T-Y2}.
This is one of the most important properties of $C_k$-move.
We consider local moves which have this
property. We define $A_1$-move as the crossing
change and $A_{k+1}$-move as a band trivialization of
$A_k$-move; see Section 1 for the precise definition.
So $A_k$-move is a generalization of $C_k$-move.
In fact, the results for $A_k$-move in this paper hold
for $C_k$-move.

Let {\em $A_k$-equivalence} denote an equivalence relation given by
$A_k$-moves and ambient isotopy.
Habiro showed that two oriented knots are $C_k$-equivalent
if and only if they have the same Vassiliev invariants of order
$\leq k-1$ \cite{Habiro1},\cite{Habiro}.
The \lq only if' part of this result is true for $A_k$-move
and for the embeddings of a graph, in particular for links
(Theorem 5.1).
However the \lq if' part does not hold for two-component links.
For example, the Whitehead link is not $C_3$-equivalent to a trivial link
because they have different Arf invariants, see \cite{T-Y}.
On the other hand, H. Murakami showed in
\cite{Murakami} that the Vassiliev invariants of links of order
$\leq 2$ are determined by the linking numbers and the
second coefficient of the Conway polynomial of each component.
Hence, the values of any Vassiliev invariant of order $\leq 2$
of these two links are the same.
So we note that Vassiliev invariants of order $\leq k-1$ are not enough to
characterize $C_k$-equivalent embeddings of a graph.

We will define in Section 1 a {\em finite type invariant of order $(n;k)$}
as a generalization of a Vassiliev invariant and
see that if any finite type invariants of order $(1;k)$ takes
the same value on two $A_{k-1}$-equivalent embeddings of a graph,
then they are $A_k$-equivalent (Theorem 1.1).
While a Vassiliev invariant is defined by the change in its value
at every \lq wall' corresponding to a crossing change,
a finite type invariant of order $(n;k)$
is defined similarly by \lq walls' corresponding to $A_k$-moves.
A finite type invariant of order $(n;1)$ is a Vassiliev
invariant of order $\leq n$.

It is shown that the set of $C_k$-equivalence classes of knots
forms an abelian group under the connected sum \cite{Habiro1},\cite{Habiro}.
This is also true for $A_k$-equivalence classes.
Since the connected sum is peculiar to knots, we cannot apply
it to embeddings of a graph. In Section 2, we will define a certain
geometric sum for the elements in an $A_{k-1}$-equivalence class
of the embedding of a graph.
Then we will see that the quotient set of the $A_{k-1}$-equivalence
class under $A_k$-equivalence forms an abelian group (Theorem 2.4).

It is not essential that $A_1$-move is the crossing change.
This is a big difference between $A_k$-move and $C_k$-move.
We will study a generalization of $A_k$-move in Section 4.
For example, if we put $A_1$-move to be the $\#$-move defined
by Murakami \cite{Murakami1}, then we get several results
similar to that for original $A_k$-move.

\bigskip
\noindent
{\bf 1. $A_k$-Moves and Finite Type Invariants}

\bigskip

Let $B^3$ be the oriented unit 3-ball.
A {\em tangle} is a
disjoint union of properly embedded arcs in $B^3$. A tangle is {\em
trivial} if it is
contained in a properly embedded 2-disk in $B^3$.
A {\em trivialization} of a tangle $T=t_1\cup t_2\cup\cdots\cup t_k$ is a
choice of
mutually disjoint disks $D_1,D_2,...,D_k$ in $B^3$ such that $\d
D_i=(D_i\cap\partial B^3)\cup t_i$ for $i=1,2,...,k$.
It can be shown that in general a trivialization is not
unique up to ambient isotopy of $B^3$ fixed on the tangle.

Let $T$ and $S$ be tangles, and let $t_1,t_2,...,t_k$ and
$s_1,s_2,...,s_k$ be the components of $T$ and $S$ respectively.
Suppose that for each
$t_i$ there exists some $s_j$ such that $\partial t_i=\partial s_j$.
Then we call the ordered pair
$(T,S)$ a {\em local move}, which can be interpreted as substituting
$S$ for $T$. Two local moves $(T,S)$ and
$(T',S')$ are {\em equivalent} if there exists an orientation preserving
homeomorphism
$h:B^3\longrightarrow B^3$ such that $h(T)=T'$ and $h(S)$
is ambient isotopic to $S'$ relative to $\partial B^3$.
We consider local moves up to this equivalence.

Let $(T,S)$ be a local move such that $T$ and $S$ are trivial tangles.
First choose a trivialization $D_1,D_2,...,D_k$ of $T$.
Each $D_i$ intersects $\partial B^3$ in an arc $\gamma_i$.
Let $E_i$ be a small regular neighbourhood of
$\gamma_i$ in $\partial B^3$. We devide the circle $\partial E_i$
into two arcs $\alpha_i$ and $\beta_i$ such that $\alpha_i\cap\beta_i=
\partial\alpha_i=\partial\beta_i$.
By slightly perturbing ${\mathrm{int}}\alpha_i$ and
${\mathrm{int}}\beta_i$ into the interior of
$B^3$ on either side of $D_i$,
we obtain properly embedded arcs $\tilde{\alpha}_i$ and $\tilde{\beta}_i$.
We consider $k$ local moves $(S\cup\tilde{\alpha}_i,S\cup\tilde{\beta}_i)\
(i=1,2,...,k)$
and call them the {\em band trivializations} of the local move $(T,S)$ with
respect to the trivialization $D_1,D_2,...,D_k$.
Note that both $S\cup\tilde{\alpha}_i$ and
$S\cup\tilde{\beta}_i$ are trivial tangles.

We now inductively define a sequence of local moves on trivial tangles
in $B^3$ which depend on the choice of trivialization.
An {\em $A_1$-move} is the crossing change shown in Figure 1.1.
Suppose that $A_k$-moves are defined and there are $l$ $A_k$-moves
$(T_1,S_1),(T_2,S_2),...,(T_l,S_l)$ up to equivalence.
For each $A_k$-move $(T_i,S_i)\ (i=1,2,...,l)$, we choose a single
trivialization $\tau_i=\{D_{i,1},D_{i,2},...,D_{i,k+1}\}$ of $T_i$ and
fix it. (The choice of $\tau_i$ is independent of the trivialization
that is chosen to define $A_k$-move $(T_i,S_i)$.)
Then the band trivializations of $(T_i,S_i)$ with respect to
the trivialization $\tau_i$ are called {\em $A_{k+1}(\tau_i)$-moves} and
these $A_{k+1}(\tau_i)$-moves $(i=1,2,...,l)$ are called
$A_{k+1}(\tau_1,\tau_2,...,\tau_l)$-moves.
Note that the number of $A_{k+1}(\tau_1,\tau_2,...,\tau_l)$-moves is
at most $l(k+1)$ up to equivalence.
Although the choice of trivializations is important for
the definition of $A_k$-move, our proof is the same for every choice.
Therefore the results of this paper hold for every choice of
trivializations $\tau_1,\tau_2,...,\tau_l$.
So we denote $A_{k+1}(\tau_1,\tau_2,...,\tau_l)$-move simply
as {\em $A_{k+1}$-move}.
It is known that $C_k$-move defined by Habiro is a
special case of $A_k$-move for certain choices of
trivializations; see \cite{Habiro2}, \cite{O-T-Y}.
We will see that $A_k$-move, as well as $C_k$-move, has the
property mentioned in Introduction (Proposition 2.1 and Lemma 2.2).

\figure{25}{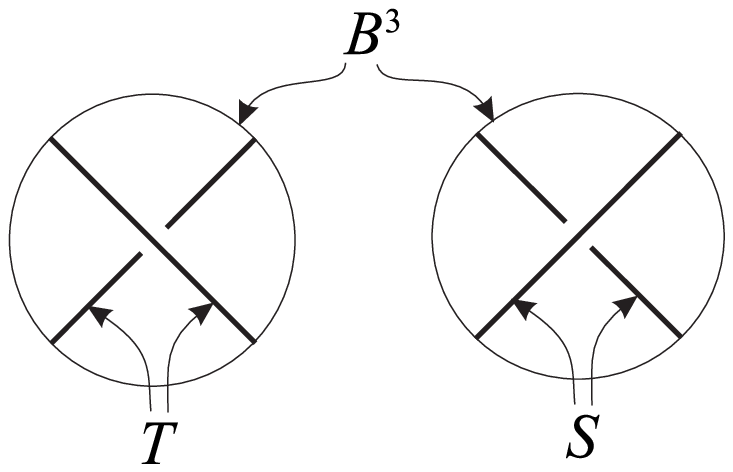}{1.1}

\medskip
\noindent
{\em Examples.} (1) The trivialization of a tangle in Figure 1.1 is
unique up to ambient isotopy.
Therefore we have any band trivialization of an
$A_1$-move is equivalent  to the local move in Figure 1.2-(i). Thus
$A_2$-move is unique up  to
equivalence. It is not hard to see that an $A_2$-move is equivalent to the
{\em delta move} in Figure 1.2-(ii) defined by H. Murakami and Y. Nakanishi
\cite{M-N}, and then it is equivalent to the local move in
Figure 1.2-(iii).

\noindent
(2) If we choose a trivialization for the $A_2$-move as in Figure 1.3-(i),
then, by the symmetry of the $A_2$-move, any $A_3$-move is
equivalent to the local move in Figure 1.3-(ii).

\begin{center}
\begin{tabular}{c}
\begin{tabular}{ccc}
\includegraphics[trim=0mm 0mm 0mm 0mm, width=.25\linewidth]
{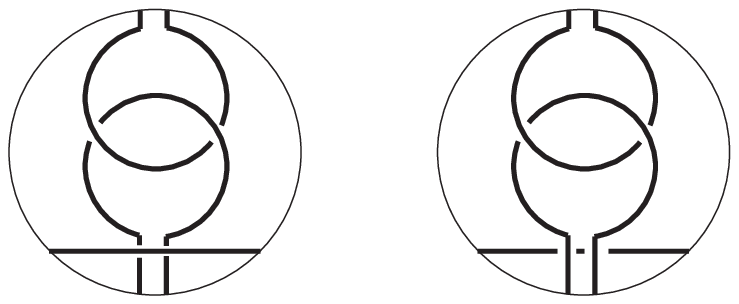}\hspace*{1cm} &
\includegraphics[trim=0mm 0mm 0mm 0mm, width=.25\linewidth]
{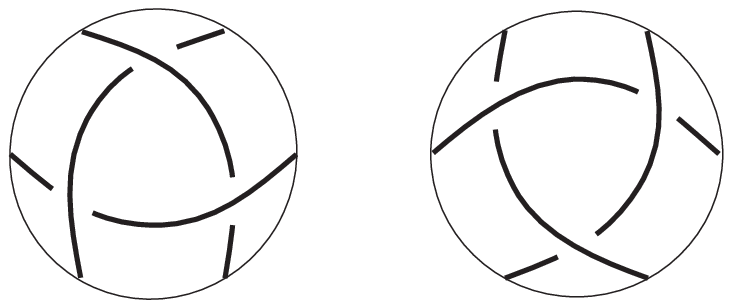}\hspace*{1cm} &
\includegraphics[trim=0mm 0mm 0mm 0mm, width=.25\linewidth]
{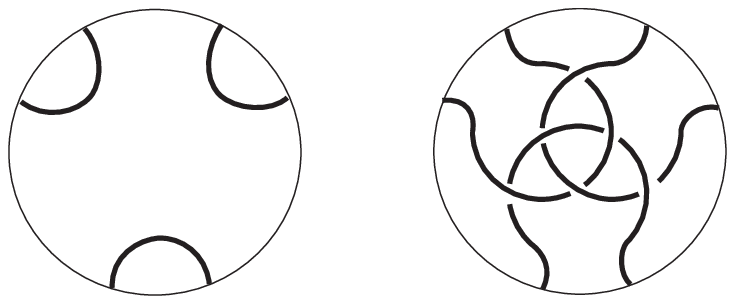}\\
(i)\hspace*{1cm} & (ii)\hspace*{1cm} & (iii)
\end{tabular}\\
Figure 1.2
\end{tabular}
\end{center}
\begin{center}
\begin{tabular}{c}
\begin{tabular}{ccc}
\includegraphics[trim=0mm 0mm 0mm 0mm, width=.25\linewidth]
{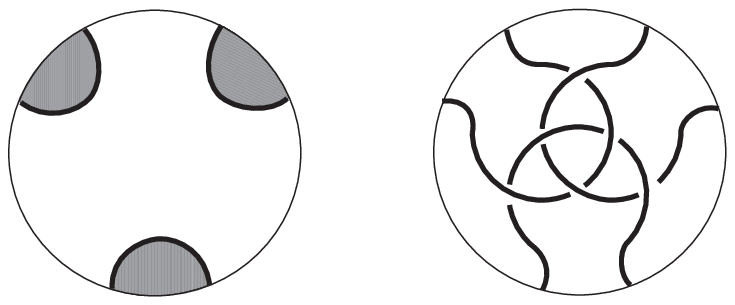}\hspace*{1cm} &
\includegraphics[trim=0mm 0mm 0mm 0mm, width=.25\linewidth]
{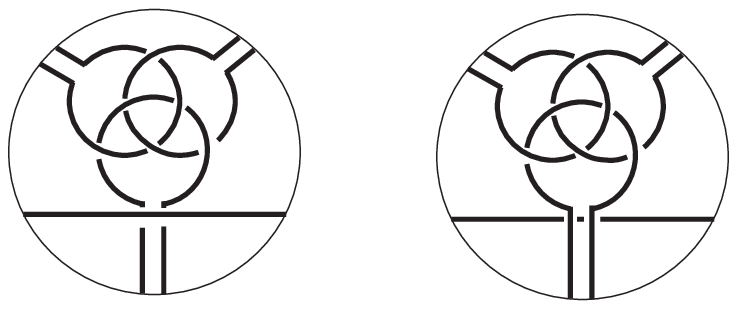}\\
(i)\hspace*{1cm} & (ii)
\end{tabular}\\
Figure 1.3
\end{tabular}
\end{center}

A local move $(S,T)$ is called the {\em inverse} of a local move
$(T,S)$. It is clear that the inverse of an $A_1$-move is again an
$A_1$-move. By the definition of $A_k$-move,
we see that the inverse of an $A_k$-move with $k\geq2$  is equivalent to
itself.

Let $(T,S)$ be an $A_k$-move and $D_1,D_2,\cdots,D_{k+1}$ the fixed
trivialization of $T=t_1\cup
t_2\cup\cdots\cup t_{k+1}$. We set $\alpha=\partial B^3\cap(D_1\cup
D_2\cup\cdots\cup D_{k+1})$ and $\beta=S$.  A link
$L$ in $S^3$ is called {\em type $k$} if there is an orientation preserving
embedding
$\varphi:B^3\longrightarrow S^3$ such that $L=\varphi(\alpha\cup\beta)$.
Then the pair $(\alpha,\beta)$ is called a {\em link model of $L$}.

We now define an equivalence relation on spatial graphs by $A_k$-move.
Let $G$ be a finite graph. Let $V(G)$ denote the set of the vertices of $G$.
Let $f,g:G\longrightarrow S^3$ be embeddings.
We say  that $f$ and $g$ are {\em related by an
$A_k$-move} if there is an
$A_k$-move $(T,S)$ and an orientation preserving embedding
$\varphi:B^3\longrightarrow S^3$ such that
\begin{enumerate}
\item[(i)] if $f(x)\neq g(x)$ then both $f(x)$ and $g(x)$ are contained in
$\varphi({\mathrm{int}} B^3)$,
\item[(ii)] $f(V(G))=g(V(G))$ is disjoint from $\varphi(B^3)$, and
\item[(iii)] $f(G)\cap\varphi(B^3)=\varphi(T)$ and $g(G)\cap\varphi(B^3)=
\varphi(S)$.
\end{enumerate}
We also say that $g$ is obtained from $f$ by an {\em application} of
$(T,S)$. We define {\em$A_k$-equivalence} as an equivalence relation
on the set of all embeddings of $G$ into $S^3$
given by the relation above and  ambient isotopy.
For an embedding $f:G\longrightarrow S^3$, let $[f]_k$ denote
the $A_k$-equivalence class of $f$.
By the definition of $A_k$-move we see that an
application of an $A_{k+1}$-move is
realized by two applications of $A_k$-move and ambient isotopy.
Thus $A_{k+1}$-equivalence implies $A_k$-equivalence.
In other words we have
$[f]_1\supset[f]_2\supset\cdots\supset[f]_k\supset[f]_{k+1}\supset\cdots$.

Let $f:G\longrightarrow S^3$ be an embedding, $L_i$ links of type
$k$ and  $(\alpha_i,\beta_i)$ their link models $(i=1,2,...,n)$.
Let $I=[0,1]$ be the unit closed interval. An embedding
$g:G\longrightarrow S^3$ is called a {\em band sum of $f$ with
$L_1,L_2,...,L_n$} if there are mutually disjoint embeddings
$b_{ij}:I\times I\longrightarrow S^3\
(i=1,2,...,n,\ j=1,2,...,k+1)$ and mutually disjoint orientation preserving
embeddings
$\varphi_i:B^3\longrightarrow S^3-f(G)$ with
$L_i=\varphi_i(\alpha_i\cup\beta_i)$ $(i=1,2,...,n)$
such that the following (i) and (ii) hold:
\begin{enumerate}
\item[(i)] $b_{ij}(I\times I)\cap f(G)=b_{ij}(I\times I)\cap
f(G-V(G))=b_{ij}(I\times \{0\})$ and
$b_{ij}(I\times I)\cap (\bigcup_{l}\varphi_l(B^3))=b_{ij}(I\times \{1\})$ is
a component of $\varphi_i(\alpha_i)$ for any
$i,j\ (i=1,2,...,n,\ j=1,2,...,k+1)$.
\item[(ii)] $f(x)=g(x)$ if $f(x)$ is not contained in
$\bigcup_{i,j}b_{ij}(I\times\{0\})$ and

$\displaystyle g(G)=(f(G)\cup\bigcup_{i}L_i
-\bigcup_{i,j}b_{ij}(I\times\partial I))
\cup\bigcup_{i,j} b_{ij}(\partial I\times I)$.
\end{enumerate}
Then we denote $g$ by $F(f;\{L_1,L_2,...,L_n\},\{B_1,B_2,...,B_n\})$,
where $B_i=b_{i1}(I\times I)\cup b_{i2}(I\times I)\cup\cdots\cup
b_{ik+1}(I\times I)\ (i=1,2,...,n)$.
We call each
$b_{ij}(I\times I)$ a {\em band}.  We call each $\varphi_i(B^3)$ an {\em
associated ball} of
$L_i$. See Figure 1.4 for an example of a band sum of an embedding
$f$ with links $L_1,L_2,L_3$ of type 3.

\figure{99}{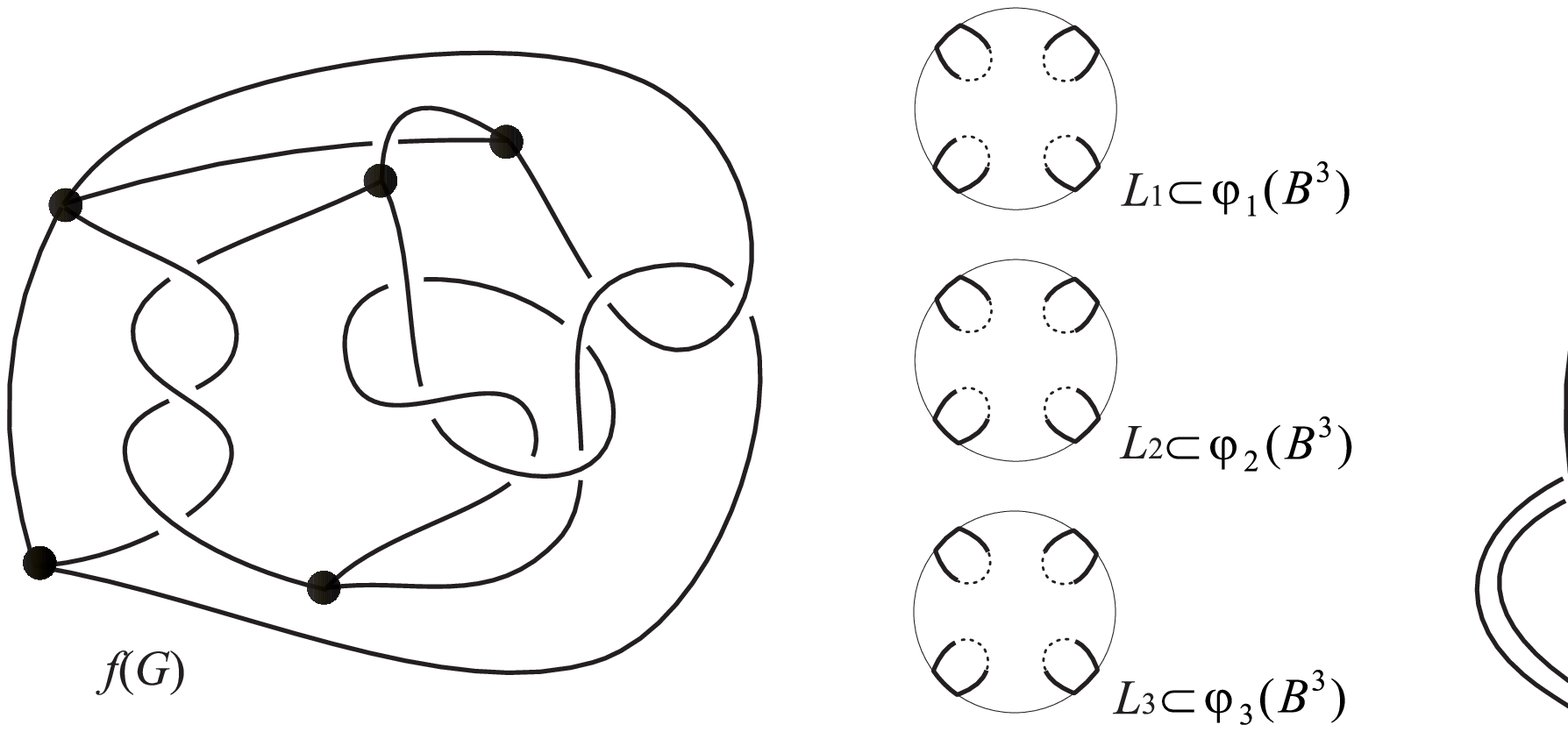}{1.4}

\medskip
\noindent
{\em Remark.} It follows from the definition that
if $g$ is a band sum of $f$ with some links of type $k$, then
$g$ is $A_k$-equivalent to $f$. The converse is also true and will be shown
in Proposition 2.1. In Lemma 2.2, we show that the position of a band is
changeable up to $A_{k+1}$-equivalence. The origin of the name
\lq band trivialization' comes from this fact.

\medskip
 Let $h:G\longrightarrow S^3$ be an embedding and $H$ an abelian
group. Let $\varphi:[h]_{k-1}\longrightarrow H$ be an invariant.
We say that $\varphi$ is a {\em finite type invariant of order $(n;k)$}
if for any embedding $f\in[h]_{k-1}$
and any band sum $F(f;\{L_1,L_2,...,L_{n+1}\},\{B_1,B_2,...,B_{n+1}\})$
of $f$ with links $L_1,L_2,...,L_{n+1}$ of type $k-1$,
\[\sum_{X\subset\{1,2,...,n+1\}}(-1)^{|X|}\varphi
(F(f;\bigcup_{i\in X}\{L_i\},
\bigcup_{i\in X}\{B_i\}))=0\in H,\]
where the sum is taken over all subsets, including the empty set, and
$|X|$ is the number of the elements in $X$.

In the next section we  show the following theorem.

\medskip\noindent
{\bf Theorem 1.1.} {\em Let $f,g:G\longrightarrow S^3$ be
$A_{k-1}$-equivalent embeddings.
Then they are $A_k$-equivalent if and only if $\varphi(f)=
\varphi(g)$ for any finite
type $A_k$-equivalence invariant $\varphi$ of order $(1;k)$.}

\medskip
Note that finite type invariants of order $(n;2)$ coincide with Vassiliev
invariants of order $n$.
It is shown in \cite[Theorem 1.1, Theorem 1.3]{M-T} that two embeddings
of a finite graph $G$ into $S^3$ are $A_2$-equivalent if and only if
they have the same Wu invariant \cite{Wu}. It follows from
\cite[Section 2]{Tani} that Wu invariant is a finite type invariant of
order $(1;2)$. Since two embeddings are always $A_1$-equivalent, we have
the following corollary.

\medskip\noindent
{\bf Corollary 1.2.} {\em Let $f,g:G\longrightarrow S^3$
be embeddings. Then the following
conditions are mutually equivalent.
\begin{enumerate}
\item[$\mathrm{(i)}$] $f$ and $g$ are $A_2$-equivalent.
\item[$\mathrm{(ii)}$] $f$ and $g$ have the same Wu invariant.
\item[$\mathrm{(iii)}$] $\varphi(f)=\varphi(g)$
for any Vassiliev invariant $\varphi$
of order $1$. \fin
\end{enumerate}}

\medskip
In Section 5 we show the following proposition.

\medskip\noindent
{\bf Proposition 1.3.} {\em Let $\varphi$ be a Vassiliev invariant of order
$(n+1)(k-1)-1$. Then $\varphi$ is a finite type invariant of order $(n;k)$.}

\bigskip
\noindent
{\bf 2. $A_k$-Equivalence Group of Spatial Graphs}

\bigskip
The following proposition is a natural generalization of
\cite[Lemma]{Yasuhara2}
and stems from the fact that a knot with the unknotting number $u$
can be unknotted by changing $u$ crossings of a regular diagram
of it \cite{Suzuki}, \cite{Yamamoto}.

\medskip\noindent
{\bf Proposition 2.1.} {\em Let $f,g:G\longrightarrow S^3$ be embeddings.
If $f$ and $g$ are $A_k$-equivalent, then $g$ is ambient isotopic to
a band sum of $f$ with some links of type $k$. }

\medskip\noindent
{\em Proof.} We consider the embeddings up to ambient isotopy for simplicity.
By the assumption there is a finite sequence of embeddings
$f=f_0,f_1,...,f_n=g$
and orientation preserving embeddings $\varphi_1,\varphi_2,...,\varphi_n:
B^3\longrightarrow S^3$ such
that $(\varphi_i^{-1}(f_{i-1}(G)),\varphi_i^{-1}(f_{i}(G)))$ is an
$A_k$-move for each $i$.  We shall prove this proposition by
induction on $n$.

First we consider the case $n=1$. Let $D_1,D_2,...,D_{k+1}$ be the fixed
trivialization
of the tangle $\varphi_1^{-1}(f_{0} (G))$ and $\gamma_j=D_j\cap\partial B^3\
(j=1,2,...,k+1)$. Then $L=\bigcup_{j}\varphi_1(\gamma_j)\cup(\varphi_1(B^3)
\cap f_1(G))$ is a link of type $k$.
By taking a small one-sided collar for each $\varphi_1(\gamma_j)$
in $S^3-\varphi_1({\mathrm{int}}B^3)$, we have
mutually disjoint embeddings $b_j:I\times I\longrightarrow
S^3\ (j=1,2,...,k+1)$ such that
$b_j(I\times I)\cap \varphi_1(B^3)=b_j(I\times\{1\})=\varphi_1(\gamma_j)$ and
$b_j(I\times I)\cap f_0(G)=b_j(I\times I)\cap f_1(G)=b_j(\partial I\times I)$.
Then we deform $f_0$ up to
ambient isotopy along the disk $b_j(I\times I)\cup\varphi_1(D_j)$ such that
$b_j(I\times I)\cap
f_0(G)=b_j(I\times\{0\})$ for each $j$. Then we have a required band sum
$g=F(f_0;\{L\},\{B\})$, where $B=b_{1}(I\times I)\cup b_{2}(I\times
I)\cup\cdots\cup b_{k+1}(I\times
I)$.

Next suppose that $n>1$. By the hypothesis of our induction, $g$ is a band
sum $F(f_1;{\cal L},{\cal B})$, where ${\cal L}=\{L_1,L_2,...,L_{n-1}\}$ is
a set of links of type
$k$, ${\cal B}=\{B_1,B_2,...,B_{n-1}\}$ and each $B_i$ is a union of bands
attaching to $L_i$. Deform
$F(f_1;{\cal L},{\cal B})$ up to ambient isotopy keeping the image $f_1(G)$
so that neither the
associated balls of
$\cal L$ nor  the bands in
$\cal B$ intersect
$\varphi_1(B^3)$.  Note that this deformation is possible,
since the tangle $\varphi_1^{-1}(f_{1}(G))$ is trivial. In fact, sweeping out
the associated balls,
band-slidings and sweeping out the bands are sufficient. See Figure 2.1.
Then by the same arguments as
that in the case $n=1$, we find that
$f_1$ is a band sum $F(f;\{L\},\{B\})$. Then we have
\[F(F(f;\{L\},\{B\});{\cal L},{\cal
B})= F(f;\{L\}\cup{\cal L},\{B\}\cup{\cal B}).\]  This completes the proof.
\fin

\figure{7}{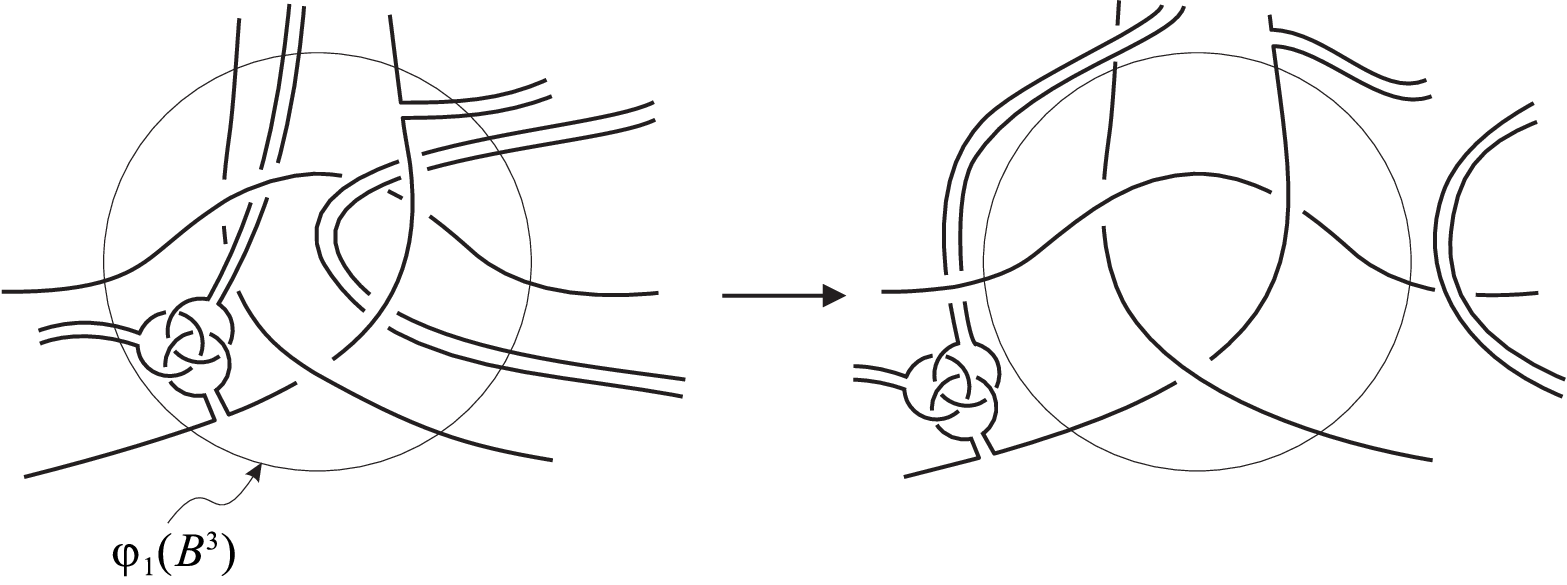}{2.1}

As we mentioned before, the origin of the name \lq band trivialization'
comes from the following lemma.

\medskip
\noindent
{\bf Lemma 2.2.}
{\em The moves in Figures $2.2$-$(\mathrm i)$, $(\mathrm ii)$,
$(\mathrm iii)$ and $(\mathrm iv)$
 are realized by $A_{k+1}$-moves. }

\begin{center}
\begin{tabular}{c}
\begin{tabular}{cc}
\includegraphics[trim=0mm 0mm 0mm 0mm, width=.31\linewidth]
{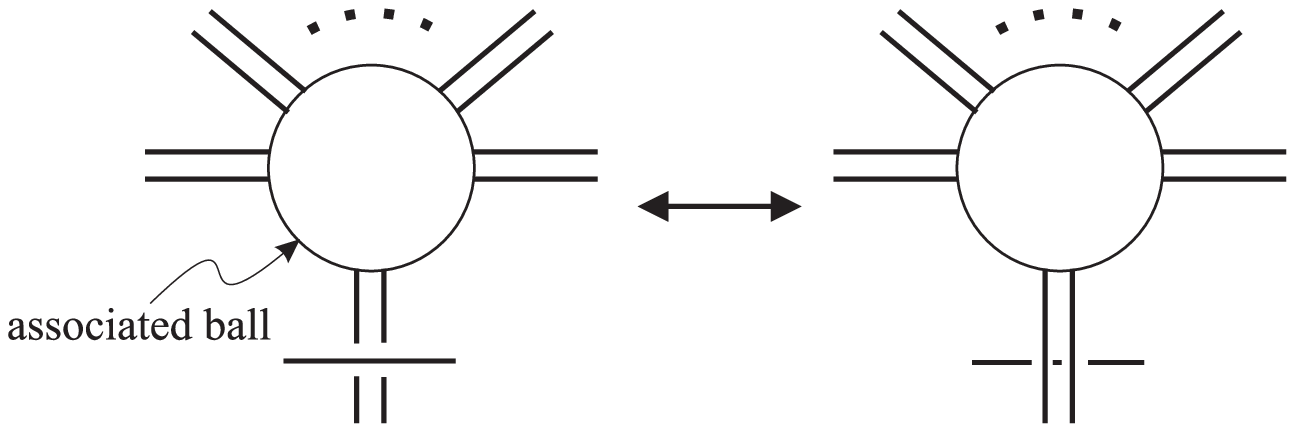} &
\includegraphics[trim=0mm 0mm 0mm 0mm, width=.276\linewidth]
{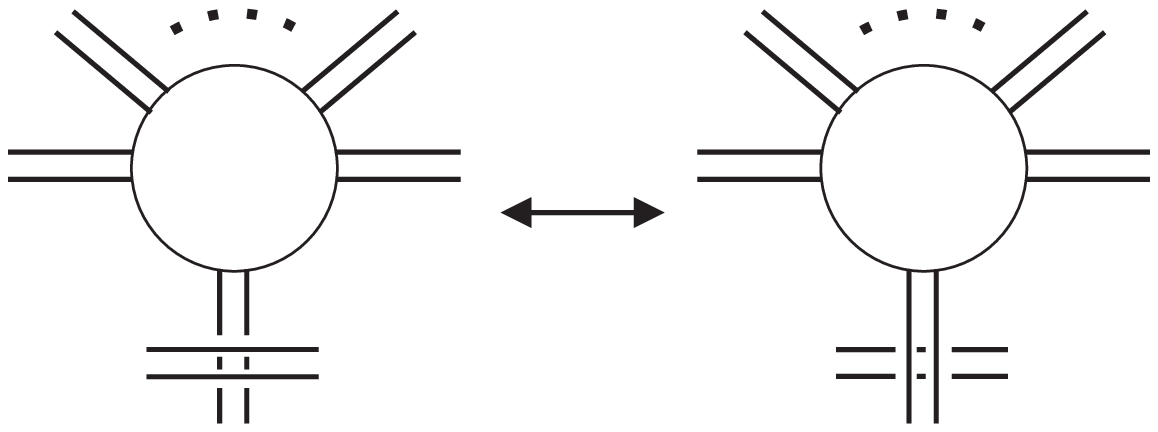}\\
\hspace*{5mm}(i) & (ii)\\
\hspace*{5mm}
\includegraphics[trim=0mm 0mm 0mm 0mm, width=.276\linewidth]
{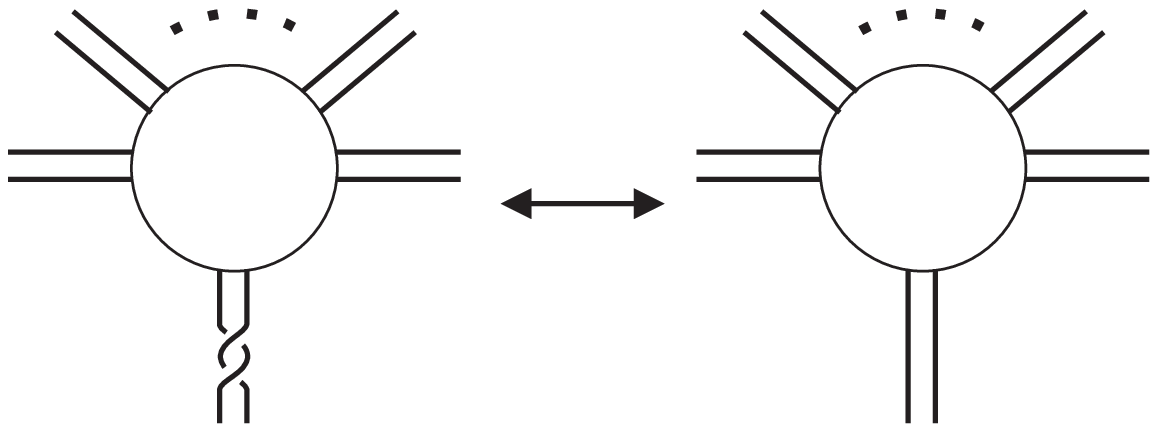} &
\includegraphics[trim=0mm 0mm 0mm 0mm, width=.5\linewidth]
{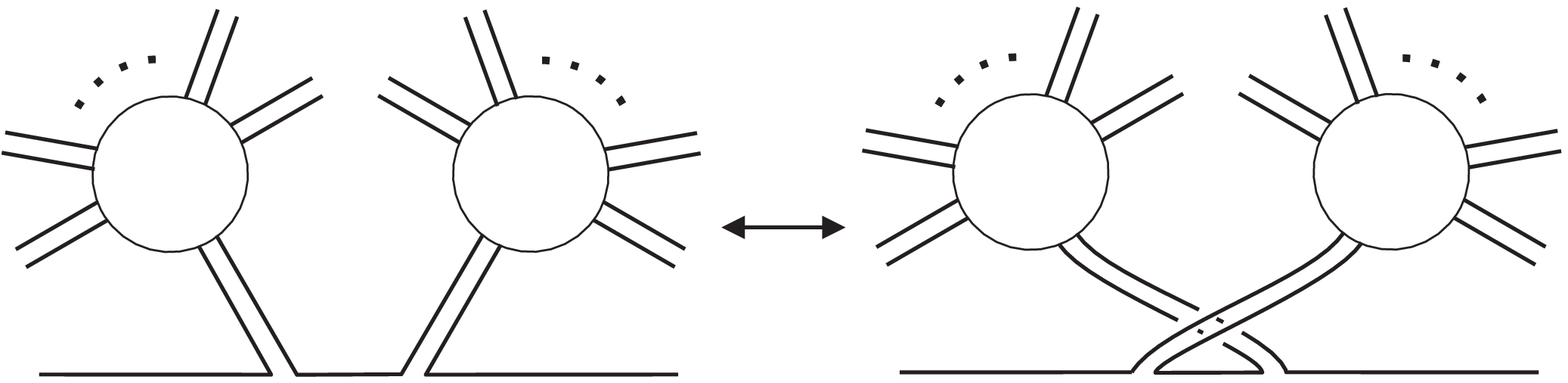}\\
\hspace*{5mm}(iii) & (iv)
\end{tabular}\\
Figure 2.2
\end{tabular}
\end{center}

\medskip
{\bf Proof.}
The move in Figure 2.2-(i) is just a band trivialization of an
$A_k$-move.  Hence by the definition it is an $A_{k+1}$-move.
It is easy to see that the moves in Figures 2.2-(ii) and (iii) are
generated by the moves in Figure 2.2-(i).  To see that the move in
Figure 2.2-(iv) is realized by $A_{k+1}$-moves,
we first slide the bands as illustrated in Figure 2.3,
and then perform the moves in Figure 2.2-(i).
\fin

\figure{6}{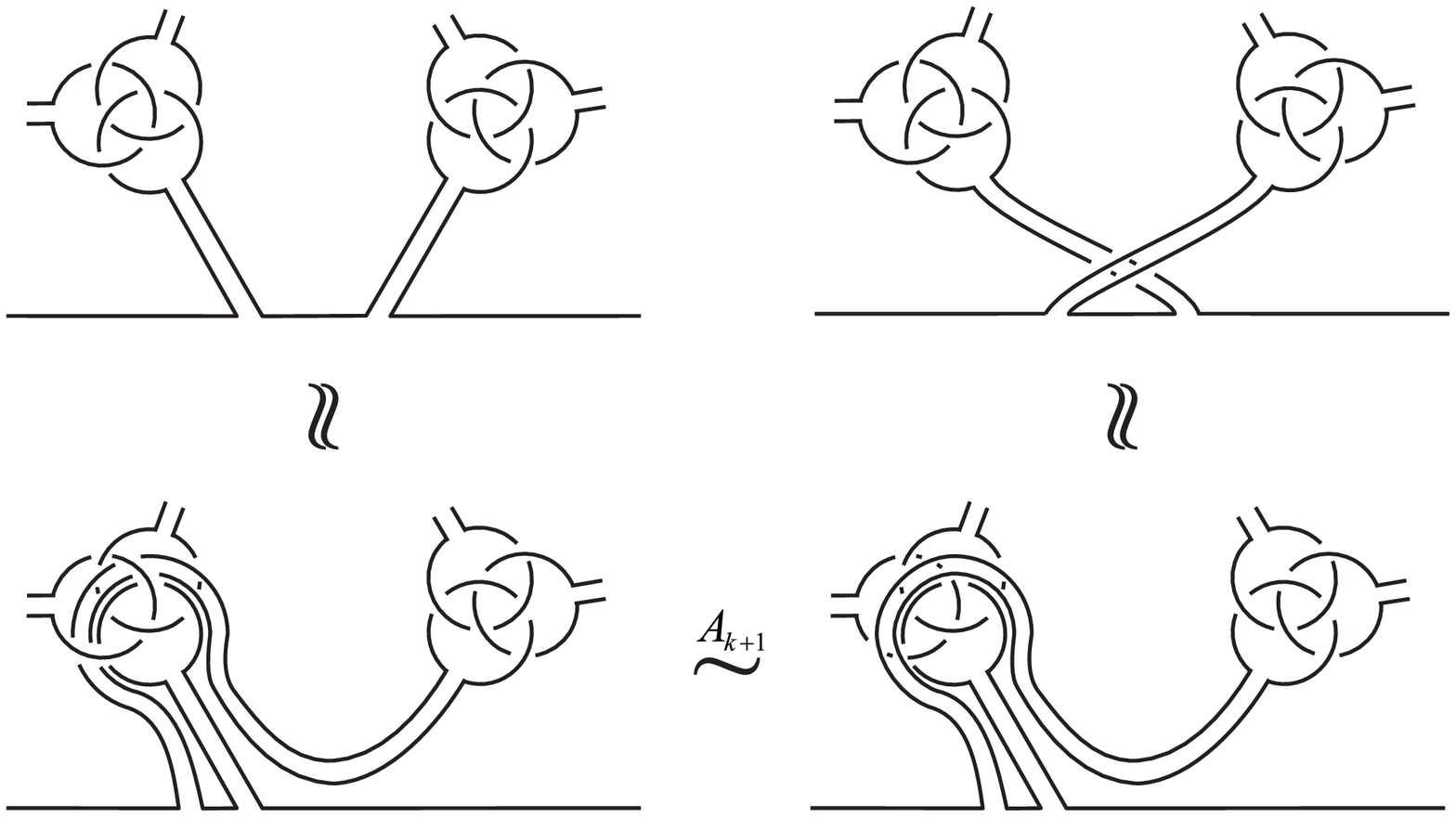}{2.3}

Let $h:G\longrightarrow S^3$ be an embedding and let
$[f_1]_k,[f_2]_k\in[h]_{k-1}/(A_k\mbox{-equivalence})$,
where $[h]_{k-1}/(A_k\mbox{-equivalence})$ denotes the set of
$A_k$-equivalence classes in $[h]_{k-1}$.
Since both $f_1$ and $f_2$ are $A_{k-1}$-equivalent to $h$, by
Proposition 2.1, there are band sums $F(h;{\cal L}_i,{\cal B}_i)\in
[f_i]_k$  of $h$ with links ${\cal L}_i$ of type $k-1$$(i=1,2)$.
Suppose that the bands in ${\cal B}_1$ and
the associated balls of ${\cal L}_1$ are disjoint from the bands in ${\cal
B}_2$ and the
associated balls of ${\cal L}_2$. Note that up to slight ambient isotopy of
$F(h;{\cal L}_2,{\cal
B}_2)$ that preserves $h(G)$ we can always choose the bands and the
associated balls so that they
satisfy this condition. In the following we assume this condition without
explicit mention. Then we
have a new band sum
$F(h;{\cal L}_1\cup{\cal L}_2,{\cal B}_1\cup{\cal B}_2)$.  We define
\[[f_1]_k+_h[f_2]_k=[F(h;{\cal L}_1\cup{\cal L}_2,{\cal B}_1\cup{\cal
B}_2)]_k.\]

\medskip\noindent
{\bf Lemma 2.3.} {\em The sum \lq\ $+_h$' above is well-defined.}

\medskip\noindent
{\em Proof.} It is sufficient to show for two embeddings $F(h;{\cal
L}_1,{\cal B}_1),
F(h;{\cal L}'_1,{\cal B}'_1)\in [f_1]_k$ that
$F(h;{\cal L}_1\cup{\cal L}_2,{\cal B}_1\cup{\cal B}_2)$ and
$F(h;{\cal L}'_1\cup{\cal L}_2,{\cal B}'_1\cup{\cal B}_2)$ are
$A_k$-equivalent. Consider a sequence of ambient isotopies and
applications of $A_k$-moves that
deforms $F(h;{\cal L}_1,{\cal B}_1)$ into $F(h;{\cal L}'_1,{\cal B}'_1)$.
We consider this sequence
of deformations together with the links in ${\cal L}_2$ and the bands in
${\cal B}_2$. Whenever we
apply an $A_k$-move we deform the associated balls of ${\cal L}_2$ and the
bands in ${\cal B}_2$ up to
ambient isotopy so that they are away from the 3-ball within which the
$A_k$-move  is applied. Thus
$F(h;{\cal L}_1\cup{\cal L}_2,{\cal B}_1\cup{\cal
B}_2)=F(F(h;{\cal L}_1,{\cal B}_1);{\cal
L}_2,{\cal B}_2)$ is $A_k$-equivalent to a band sum $F(F(h;{\cal
L}'_1,{\cal B}'_1);{\cal
L}'_2,{\cal B}'_2)$ for some ${\cal L}'_2$ and ${\cal B}'_2$. Compare the
band sums $F(F(h;{\cal
L}'_1,{\cal B}'_1);{\cal L}'_2,{\cal B}'_2)$ and $F(h;{\cal L}'_1\cup{\cal
L}_2,{\cal B}'_1\cup{\cal
B}_2)=F(F(h;{\cal L}'_1,{\cal B}'_1);{\cal L}_2,{\cal B}_2)$. We have
that the links in ${\cal L}'_2$
are ambient isotopic to the links in ${\cal L}_2$.
It follows from Lemma 2.2 that the bands in ${\cal B}'_2$
can be deformed into the position of the bands in ${\cal B}_2$
by band slidings and $A_k$-moves.
Thus these two are $A_k$-equivalent. \fin

\medskip\noindent
{\bf Theorem 2.4.} {\em The set $[h]_{k-1}/(A_k\mbox{-equivalence})$ forms
an abelian group
under \lq\ $+_h$' with the unit element $[h]_k$.}

\medskip\noindent
We denote this group by ${\cal G}_k(h;G)$ and call it the
{\em $A_k$-equivalence
group} of the spatial embeddings of $G$ with the unit element $[h]_k$.

\medskip\noindent
{\em Remark.} Note that for any graph $G$ and any embedding
$h:G\longrightarrow S^3$,
$[h]_1$ is equal to the set of all embeddings of $G$ into $S^3$. In
\cite{Yasuhara}, the second
author called ${\cal G}_2(h;G)$ a {\em graph homology group} and
gave a practical method of calculating this group.

\medskip\noindent
{\em Proof.} We consider embeddings up to ambient isotopy for simplicity.
It is sufficient to show that for any $[f]_k\in
[h]_{k-1}/(A_k\mbox{-equivalence})$, there is an inverse of $[f]_k$.
Since $f$ and $h$ are
$A_{k-1}$-equivalent, by Proposition 2.1, $f$ and $h$ are band sums
$F(h;{\cal L},{\cal B})$ and $F(f;{\cal L}',{\cal B}')$ respectively,
where ${\cal L}$ and ${\cal L}'$ are sets of links of type $k-1$. Thus we
have
$h=F(F(h;{\cal L},{\cal B});{\cal L}',{\cal B}')$.
Then, by using Lemma 2.2, we deform the
associated balls of ${\cal L}'$ and the bands in ${\cal B}'$ up to
$A_k$-equivalence so that they
are disjoint form the associated balls of ${\cal L}$ and the bands in
${\cal B}$. Thus we see that
$h=F(F(h;{\cal L},{\cal B});{\cal L}',{\cal B}')$
is $A_k$-equivalent to a band sum
$F(h;{\cal L}\cup{\cal L}'',{\cal B}\cup{\cal B}'')$ for some ${\cal L}''$
and ${\cal B}''$ (for
example see Figure 2.4).   Thus we have
\[[f]_k+_h[F(h;{\cal L}'',{\cal B}'')]_k=
[F(h;{\cal L},{\cal B})]_k+_h[F(h;{\cal L}'',{\cal B}'')]_k=
[F(h;{\cal L}\cup{\cal L}'',{\cal B}\cup{\cal B}'')]_k=[h]_k.\]
This implies that $[F(h;{\cal L}'',{\cal B}'')]_k$ is an
inverse of $[f]_k$.  \fin

\figure{5}{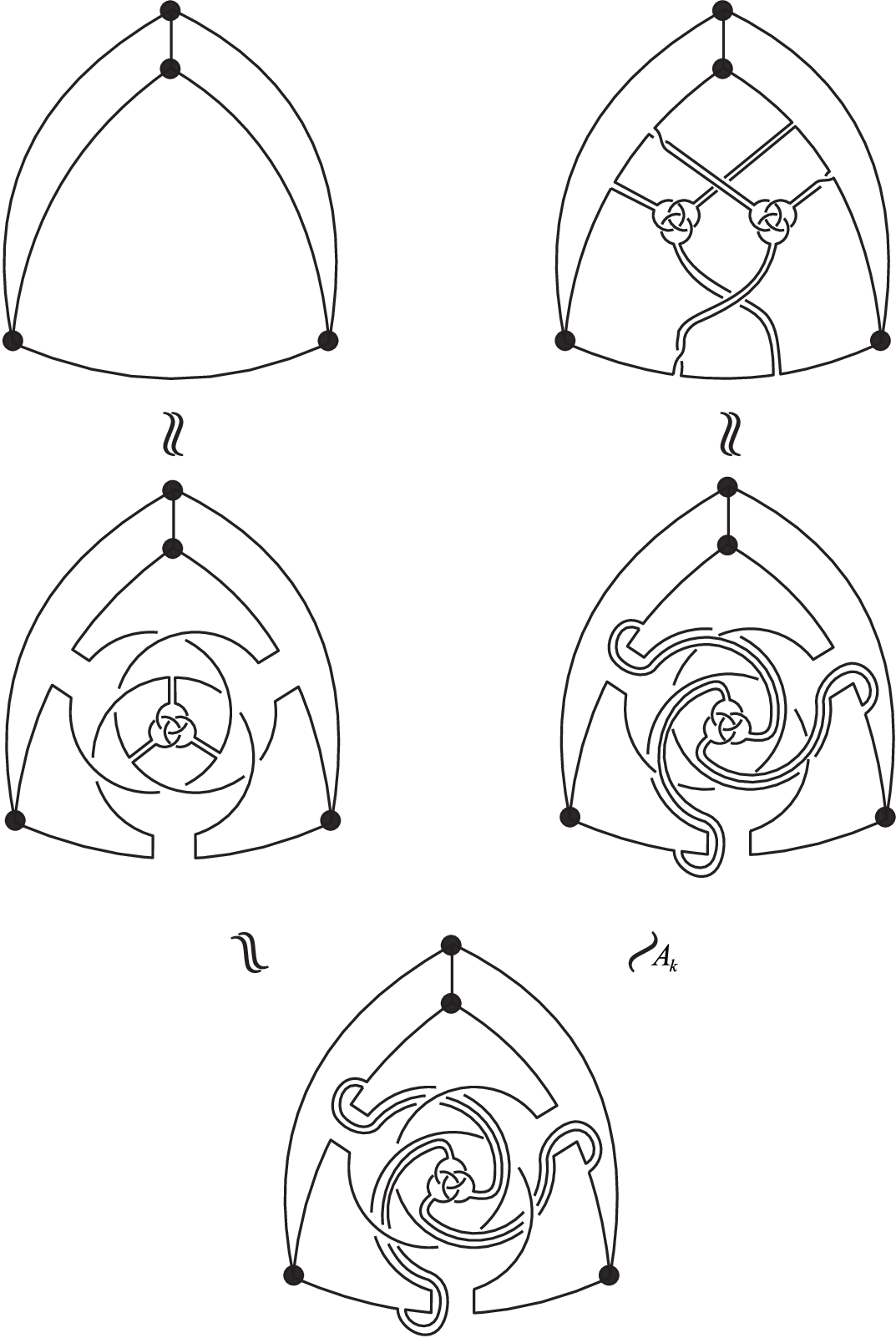}{2.4}

\medskip\noindent
{\bf Theorem 2.5.} {\em Let $h_1,h_2:G\longrightarrow S^3$ be
$A_{k-1}$-equivalent embeddings. Then the
groups ${\cal G}_k(h_1;G)$ and ${\cal G}_k(h_2;G)$ are isomorphic.}

\medskip\noindent
{\em Proof.}
We define a map $\phi:{\cal G}_k(h_1;G)\longrightarrow
{\cal G}_k(h_2;G)$ by $\phi([f]_k)=[f]_k-_{h_2}[h_1]_k$, where
$[x]_k-_{h_2}[y]_k$ denotes $[x]_k+_{h_2}(-[y]_k)$.
Clearly this map is a bijection.
We shall prove that $\phi$ is a homomorphism.
Let $[f_i]_k\in{\cal G}_k(h_1;G)\ (i=1,2)$. Then $f_i=F(h_1;{\cal
L}_i,{\cal B}_i)$ where ${\cal
L}_i$ is a set of links of type $k-1$ ($i=1,2$). Since $h_1$ and $h_2$ are
$A_{k-1}$-equivalent we see that
$h_1=F(h_2;{\cal L},{\cal B})$ where ${\cal L}$ is a set of links of type
$k-1$. Thus we have
$f_i=F(F(h_2;{\cal L},{\cal B});{\cal L}_i,{\cal B}_i)\ (i=1,2)$.
By using Lemma 2.2, we deform $f_i$ up to
$A_k$-equivalence so that the associated balls of ${\cal L}_i$ and the
bands in ${\cal B}_i$ are
disjoint from the associated balls of ${\cal L}$ and the bands in ${\cal
B}$ for $i=1,2$. We may
further assume that the associated balls of ${\cal L}_1$ and the bands in
${\cal B}_1$ are
disjoint from the associated balls of ${\cal L}_2$ and the bands in ${\cal
B}_2$. Then we have
\[\begin{array}{rl}
\phi([f_1]_k+_{h_1}[f_2]_k)&=\phi([F(F(h_2;{\cal L},{\cal B});
{\cal L}_1\cup{\cal L}_2,{\cal B}_1\cup{\cal B}_2)]_k)\\
 &=[F(F(h_2;{\cal L},{\cal B});
{\cal L}_1\cup{\cal L}_2,{\cal B}_1\cup{\cal B}_2)]_k-_{h_2}[h_1]_k\\
 &=[F(h_2;{\cal L}\cup{\cal L}_1\cup{\cal L}_2,{\cal B}\cup
{\cal B}_1\cup{\cal B}_2)]_k-_{h_2}[F(h_2;{\cal L},{\cal B})]_k\\
 &=[F(h_2;{\cal L}_1\cup{\cal L}_2,{\cal B}_1\cup
{\cal B}_2)]_k,
\end{array}\]
and for each $i\ (i=1,2)$,
\[\begin{array}{rl}
 \phi([f_i]_k)&=\phi([F(F(h_2;{\cal L},{\cal B});
{\cal L}_i,{\cal B}_i]_k)\\
 &=[F(F(h_2;{\cal L},{\cal B});{\cal L}_i,{\cal B}_i]_k-_{h_2}[h_1]_k\\
 &=[F(h_2;{\cal L}\cup{\cal L}_i,{\cal B}\cup{\cal B}_i)]_k-_{h_2}
[F(h_2;{\cal L},{\cal B})]_k\\
 &=[F(h_2;{\cal L}_i,{\cal B}_i)]_k.
\end{array}\]
Thus we have
$\phi([f_1]_k+_{h_1}[f_2]_k)=\phi([f_1]_k)+_{h_2}\phi([f_2]_k)$. \fin

\medskip\noindent
{\bf Proposition 2.6.} {\em The projection $p:[h]_{k-1}\longrightarrow
[h]_{k-1}/(A_k\mbox{-equivalence})={\cal G}_k(h;G)$ is a finite type
$A_k$-equivalence
invariant of order $(1;k)$. }

\medskip\noindent
{\em Proof.} It is clear that $p$ is an $A_k$-equivalence invariant.
We shall prove that $p$ is finite type of order $(1;k)$.
Let $f\in[h]_{k-1}$ be an embedding and $F(f;\{L_1,L_2\},\{B_1,B_2\})$
a band sum of $f$ with links $L_1,L_2$ of type $k-1$.
Then it is sufficient to show that
\[\sum_{X\subset\{1,2\}}(-1)^{|X|}p(F(f;\bigcup_{i\in
X}\{L_i\},\bigcup_{i\in X}\{B_i\}))=[h]_k.\]
Let $\phi:{\cal G}_k(f;G)\longrightarrow{\cal G}_k(h;G)$
be the isomorphism defined by
$\phi([g]_k)=[g]_k-_{h}[f]_k$. Then we have
\[\begin{array}{l}
\phi([F(f;\emptyset,\emptyset)]_k-_f[F(f;\{L_1\},\{B_1\})]_k
-_f[F(f;\{L_2\},\{B_2\})]_k+_f[F(f;\{L_1,L_2\},\{B_1,B_2\})]_k)\\
=([F(f;\emptyset,\emptyset)]_k-_h[f]_k)-_h([F(f;\{L_1\},\{B_1\})]_k-_h[f]_k)\\
\hfill-_h([F(f;\{L_2\},\{B_2\})]_k-_h[f]_k)+_h([F(f;\{L_1,L_2\},\{B_1,B_2\})
]_k-_h
[f]_k)\\
=[F(f;\emptyset,\emptyset)]_k-_h[F(f;\{L_1\},\{B_1\})]_k
-_h[F(f;\{L_2\},\{B_2\})]_k+_h[F(f;\{L_1,L_2\},\{B_1,B_2\})]_k\\
=\displaystyle{\sum_{X\subset\{1,2\}}(-1)^{|X|}p(F(f;\bigcup_{i\in
X}\{L_i\},\bigcup_{i\in X}\{B_i\}))}.
\end{array}\]
Since
\[[F(f;\{L_1,L_2\},\{B_1,B_2\})]_k=[F(f;\{L_1\},\{B_1\})]_k+_f[F(f;\{L_2\},
\{ B_2\})]_k,\]
we have
\[\begin{array}{l}
\phi([F(f;\emptyset,\emptyset)]_k-_f[F(f;\{L_1\},\{B_1\})]_k
-_f[F(f;\{L_2\},\{B_2\})]_k+_f[F(f;\{L_1,L_2\},\{B_1,B_2\})]_k)\\
=\phi([f]_k)=[h]_k.
\end{array}\]
This completes the proof. \fin

\medskip\noindent
{\em Proof of Theorem 1.1.} The \lq only if' part is clear. We show the
\lq if' part. Let $f$ and $g$ be embeddings in $[h]_{k-1}$.
Suppose that any finite type invariant of order $(1;k)$ takes the
same value on $f$ and $g$.
Then by Proposition 2.6 we have $p(f)=p(g)$, where
$p: [h]_{k-1}\longrightarrow [h]_{k-1}/(A_k\mbox{-equivalence})=
{\cal G}_k(h;G)$ is the projection.
Hence we have $[f]_k=[g]_k$. This completes the proof. \fin

\bigskip
\noindent
{\bf 3. $A_k$-Equivalence Group of Knots}

\bigskip
In this section we only consider the case that the graph $G$ is
homeomorphic to a disjoint union of circles. Let $G=S_1^1\cup
S_2^1\cup\cdots\cup S_\mu^1$.
Then there is a natural correspondence between the ambient
isotopy classes of the embeddings
of
$G$ into $S^3$ and the ambient isotopy classes of the ordered oriented
$\mu$-component links in
$S^3$. Therefore instead of specifying an embedding
$h:S_1^1\cup S_2^1\cup\cdots\cup S_\mu^1\longrightarrow S^3$,
we denote by $L$ the image
$h(S_1^1\cup S_2^1\cup\cdots\cup S_\mu^1)$ and consider it together with
the orientation of each
component and the ordering of the components. Thus ${\cal G}_k(L)$ denotes
the $A_k$-equivalence group
${\cal G}_k(h;S_1^1\cup S_2^1\cup\cdots\cup S_\mu^1)$ with the unit element
$[h]_k$.

\medskip\noindent
{\bf Theorem 3.1.} {\em Let $O$ be a trivial knot. Then for any oriented
knot $K$,
${\cal G}_k(O)$ and ${\cal G}_k(K)$ are isomorphic. }

\medskip\noindent
{\em Remark.} For a graph $G(\neq S^1)$ and embeddings
$h,h':G\longrightarrow S^3$,
${\cal G}_k(h;G)$ and  ${\cal G}_k(h';G)$ are not always isomorphic. In fact
there are two-component links $L_1$ and $L_2$ such that
${\cal G}_3(L_1)\cong {\Bbb Z}\oplus{\Bbb Z}$ and
${\cal G}_3(L_2)\cong {\Bbb Z}\oplus{\Bbb Z}\oplus{\Bbb Z}_2$ \cite{T-Y}.

\medskip\noindent
{\em Proof.} We define a map $\phi:{\cal G}_k(O)\longrightarrow
{\cal G}_k(K)$ by
$\phi([F(O;{\cal L},{\cal B})]_k)=[K\# F(O;{\cal L},{\cal B})]_k$ for each
$[F(O;{\cal L},{\cal B})]_k\in {\cal G}_k(O)$, where $\cal L$
is a set of links
of type $k-1$ and $\#$ means the connected sum of oriented knots. Clearly
this is well-defined.
By Lemma 2.2, any band sum $F(K;{\cal L},{\cal B})$ of $K$ with
links $\cal L$ of type $k-1$
is $A_k$-equivalent to  $K\#F(O;{\cal L}',{\cal B}')$ for some
links ${\cal L}'$ of type
$k-1$ and ${\cal B}'$. Hence $\phi$ is surjective.
For $[F(O;{\cal L}_i,{\cal B}_i)]_k\in{\cal G}_k(O)\ (i=1,2)$, we have
\[\begin{array}{rl}
\phi([F(O;{\cal L}_1,{\cal B}_1)]_k+_O[F(O;{\cal L}_2,{\cal B}_2)]_k)&=
\phi([F(O;{\cal L}_1\cup{\cal L}_2,{\cal B}_1\cup{\cal B}_2)]_k)\\
&=[K\#F(O;{\cal L}_1\cup{\cal L}_2,{\cal B}_1\cup{\cal B}_2)]_k\\
&=[F(K;{\cal L}_1\cup{\cal L}_2,{\cal B}_1\cup{\cal B}_2)]_k\\
&=[F(K;{\cal L}_1,{\cal B}_1)]_k+_K[F(K;{\cal L}_2,{\cal B}_2)]_k\\
&=[K\# F(O;{\cal L}_1,{\cal B}_1)]_k+_K[K\# F(O;{\cal L}_2,{\cal B}_2)]_k\\
&=\phi([F(O;{\cal L}_1,{\cal B}_1)]_k)+_K\phi([F(O;{\cal L}_2,{\cal B}_2)]_k).
\end{array}\]
This implies that $\phi$ is a homomorphism. In order to complete the proof,
we show that $\phi$ is injective.
Suppose that $[K\# F(O;{\cal L},{\cal B})]_k=[K]_k$. By  Lemma 3.2,
there is a knot $K'$  such that $[K'\# K]_k=[O]_k$. Then we have
\[[F(O;{\cal L},{\cal B})]_k=[(K'\#K)\# F(O;{\cal L},{\cal B})]_k=
[K'\#(K\# F(O;{\cal L},{\cal B}))]_k=[K'\# K]_k=[O]_k.\]
This implies that $\ker\phi=\{[O]_k\}$. \fin

\medskip
Habiro originated \lq clasper theory' and
showed  Lemma 3.2 for $C_k$-moves \cite{Habiro1},\cite{Habiro}.
The following proof is a translation  of his proof
in terms of band sum description of knots.

\medskip\noindent
{\bf Lemma 3.2.} {\em
For any knot $K$ and any integer $k\geq1$, there is a knot $K'$ such that
$K'\# K$ is $A_k$-equivalent to a trivial knot.}

\medskip\noindent
{\em Proof.} We shall prove this by induction on $k$. The case $k=1$ is
clear. Suppose that there
is a knot $K'$ such that
$K'\# K$ is $A_{k-1}$-equivalent to a trivial knot $O$ $(k>1)$. By
Proposition 2.1, we may
assume that $O=F(K'\# K;{\cal L},{\cal B})$, where $\cal L$ is a set of
links of
type $k-1$. Then, by Lemma 2.2, we see that $F(K'\# K;{\cal L},{\cal B})$
is $A_k$-equivalent to some $K \# F(K';{\cal L},{\cal B'})$.
This completes the proof. \fin

Let ${\cal K}_k$ be the set of
$A_k$-equivalence classes of all oriented knots. For $[K]_k, [K']_k\in{\cal
K}_k$, we
define $[K]_k+[K']_k=[K\# K']_k$. Then the following,
shown by Habiro \cite{Habiro1},\cite{Habiro} in the case that
$A_k$-moves coincide with $C_k$-moves,  is an immediate
consequence of Lemma 3.2.

\medskip\noindent
{\bf Theorem 3.3.} {\em
The set ${\cal K}_k$ forms an abelian group under \lq\ $+$'  with the unit
element $[O]_k$,
where $O$ is a trivial knot.  \fin  }

\bigskip
\noindent
{\bf 4. Generalized $A_k$-Move}

\bigskip
In this section, we define a generalized $A_k$-move. For this move, several
results
similar to that in Sections 1, 2 and 3 hold.

Let $T$ and $S$ be trivial tangles such that $(T,S)$ and $(S,T)$
are equivalent. Let $t_1,t_2,...,t_n$ and
$s_1,s_2,...,s_n$ be the components of $T$ and $S$ respectively.
An {\em $A_1(T,S)$-move} is this local move $(T,S)$.
Suppose that $A_k(T,S)$-moves are defined.
For each $A_k(T,S)$-move $(T_k,S_k)$ we choose a trivialization of $T_k$ and
fix it. Then the band trivializations of $(T_k,S_k)$ with respect to the
trivialization are
called {\em $A_{k+1}(T,S)$-moves}.
Let $(T_k,S_k)$ be an $A_k(T,S)$-move and $D_1,D_2,\cdots,D_{n+k-1}$
the fixed trivialization of $T_k=t_1\cup
t_2\cup\cdots\cup t_{n+k-1}$. We set $\alpha=\partial B^3\cap(D_1\cup
D_2\cup\cdots\cup D_{n+k-1})$ and
$\beta=S_k$. A link
$L$ in $S^3$ is called {\em type $(k;(T,S))$} if there is an orientation
preserving embedding $\varphi:B^3\longrightarrow S^3$
such that $L=\varphi(\alpha \cup\beta)$.
Then the pair $(\alpha,\beta)$ is called a {\em link model of $L$}.
As in Section 1,  $A_k(T,S)$-move gives an equivalence relation,
{\em $A_k(T,S)$-equivalence},  on the set of all embeddings of $G$ into
$S^3$.
For an embedding $f:G\longrightarrow S^3$, let $[f]_k$ denote
the $A_k(S,T)$-equivalence class of $f$.
 Let $h:G\longrightarrow S^3$ be an embedding and $H$ an abelian
group. Let $\varphi:[h]_{k-1}\longrightarrow H$ be an invariant.
We can define that $\varphi$ is a {\em finite type invariant of order
$(n;k;(T,S))$}
as in Section~1.

By the arguments similar to that in Sections 1, 2 and 3, we have the
following five theorems.

\medskip\noindent
{\bf Theorem 4.1.} {\em Let $f,g:G\longrightarrow S^3$ be
$A_{k-1}(T,S)$-equivalent embeddings.
Then they are $A_k(T,S)$-equivalent if and only if $\varphi(f)=\varphi(g)$
for any finite
type $A_k(T,S)$-equivalence invariant $\varphi$ of order $(1;k;(T,S))$. \fin}

\medskip
Let $h:G\longrightarrow S^3$ be an embedding.
For $[f_1]_k,[f_2]_k\in[h]_{k-1}/(A_k(T,S)\mbox{-equivalence})$,
we can define
$[f_1]_k+_h[f_2]_k$ as in Section 2, and we have

\medskip\noindent
{\bf Theorem 4.2.} {\em The set $[h]_{k-1}/(A_k(T,S)\mbox{-equivalence})$
forms an abelian group under \lq~$+_h$' with the unit element $[h]_k$. \fin}

We denote this group by ${\cal G}_{k(T,S)}(h;G)$ and call it the
{\em $A_k(T,S)$-equivalence group} of the spatial embeddings of $G$
with the unit element $[h]_k$.

\medskip\noindent
{\bf Theorem 4.3.} {\em Let $h_1,h_2:G\longrightarrow S^3$ be
$A_{k-1}(T,S)$-equivalent embeddings. Then the
groups ${\cal G}_{k(T,S)}(h_1;G)$ and ${\cal G}_{k(T,S)}(h_2;G)$ are
isomorphic. \fin}

For an embedding $h:S^1\longrightarrow S^3$, let
$K=h(S^1)$ and let ${\cal G}_{k(T,S)}(K)$ denote the $A_k(T,S)$-equivalence
group
${\cal G}_{k(T,S)}(h;S^1)$ with the unit element $[h]_k$.

\medskip\noindent
{\bf Theorem 4.4.} {\em Let $O$ be a trivial knot.
If any two knots are $A_1(T,S)$-equivalent, then for any oriented
knot $K$, ${\cal G}_{k(T,S)}(O)$ and ${\cal G}_{k(T,S)}(K)$ are isomorphic.
\fin}

Let ${\cal K}_{k(T,S)}$ be the set of
$A_k(T,S)$-equivalence classes of all oriented knots. For $[K]_k, [K']_k
\in{\cal K}_{k(T,S)}$, we define $[K]_k+[K']_k=[K\# K']_k$.

\medskip\noindent
{\bf Theorem 4.5.} {\em If any two knots are $A_1(T,S)$-equivalent,
then the set ${\cal K}_{k(T,S)}$ forms an abelian group under
\lq\ $+$'  with the unit element $[O]_k$,
where $O$ is a trivial knot.  \fin  }

\medskip\noindent
{\em Remark.} If $(T,S)$ is the $\#$-move defined by Murakami
\cite{Murakami1}, then ${\cal K}_{k(T,S)}$ is an abelian group.

\bigskip
\noindent
{\bf 5. $A_k$-Moves and Vassiliev Invariants}

\bigskip
Let $G$ be a finite graph. We give and fix orientations of the edges of
$G$. Let ${\cal
E}$ be the set of the ambient isotopy classes of the embeddings of $G$ into
$S^3$. Let ${\Bbb Z}{\cal E}$ be the free abelian group generated by
the elements of ${\cal E}$. A {\em
crossing vertex} is a double point of a map from $G$ to $S^3$ as in Figure
5.1. An {\em
$i$-singular embedding} is a map from $G$ to $S^3$ whose
multiple points are exactly
$i$ crossing vertices. By the formula in Figure 5.2 we identify an
$i$-singular embedding with an
element in ${\Bbb Z}{\cal E}$.  Let ${\cal R}_i$
be the subgroup of
${\Bbb Z}{\cal E}$ generated by all $i$-singular embeddings. Note that ${\cal
R}_i$ is independent of the
choices of the edge orientations. Let $H$ be an abelian group. Let
$\varphi:{\cal E}\longrightarrow H$ be a map.
Let $\tilde{\varphi}:{\Bbb Z}{\cal E}\longrightarrow H$
be the natural extension of
$\varphi$. We say that $\varphi$ is a {\em
Vassiliev invariant of order $n$} if $\tilde{\varphi}({\cal R}_{n+1})=\{0\}$.
Let $\iota:{\cal
E}\longrightarrow{\Bbb Z}{\cal E}$ be the natural inclusion map and
$\pi_i:{\Bbb Z}{\cal E}\longrightarrow{\Bbb Z}{\cal E}/{\cal R}_i$
the quotient homomorphism. Let $u_{i-1}=\pi_i\circ\iota:{\cal
E}\longrightarrow{\Bbb Z}{\cal E}/{\cal R}_i$ be the
composition map. Then $\varphi$ is a Vassiliev invariant of
order $n$ if and only if there is a
homomorphism $\hat{\varphi}:{\Bbb Z}{\cal E}/{\cal R}_{n+1}\longrightarrow H$
such that $\varphi=\hat{\varphi}\circ u_{n}$.
In the following we
sometimes do not distinguish between an embedding and
its ambient isotopy class so long as no confusion
occurs.

\begin{center}
\begin{tabular}{cc}
\includegraphics[trim=0mm 0mm 0mm 0mm, width=.1\linewidth]
{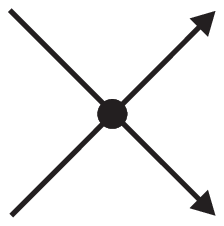}\hspace*{1cm} &
\includegraphics[trim=0mm 0mm 0mm 0mm, width=.5\linewidth]
{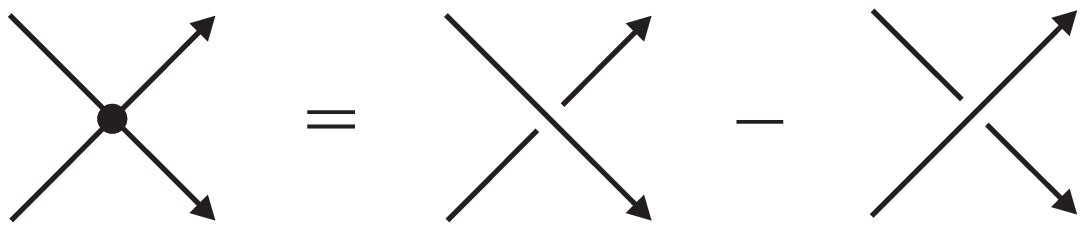}\\
Figure 5.1\hspace*{1cm} & Figure 5.2
\end{tabular}
\end{center}

\medskip\noindent
{\bf Theorem 5.1.} {\em
Let $f,g:G\longrightarrow S^3$ be $A_{k+1}$-equivalent embeddings.
Then $u_k(f)=u_k(g)$.}

\medskip
By using induction on $k$, we see that an $A_k$-move $(T,S)$ is
a $(k+1)$-component Brunnian local move, i.e., $T-t$ and $S-s$ are
ambient isotopic in $B^3$ relative $\partial B^3$ for any $t\in T$ and
$s\in S$ with $\partial t=\partial s$ \cite{T-Y0}.
It is not hard to see that if two embeddings $f$ and $g$ are related by
a $(k+1)$-component Brunnian local move, then $f$ and $g$ are
{\em $k$-similar}, where $k$-similar is
an equivalence relation defined by the first author \cite{Tani0}.
Therefore, we note that  Theorem 5.1 follows from \cite{Gusarov}
or \cite{Ohyama}. However we
give a self-contained proof here.

Let $T$ be a tangle. Let ${\cal H}(T)$ be the set of all (possibly
nontrivial) tangles that are
homotopic to $T$ relative to $\partial B^3$. Let ${\cal E}(T)$ be the
quotient of ${\cal H}(T)$ by
the ambient isotopy relative to $\partial B^3$. Then ${\Bbb Z}{\cal E}(T)$,
$i$-singular tangles and
${\cal R}_i(T)\subset{\Bbb Z}{\cal E}(T)$ are defined as above.

\medskip\noindent
{\em Proof.} It is sufficient to show for each $A_k$-move
$(T,S)$ that $T-S$ is an
element in ${\cal R}_k(T)$. We show this by induction on $k$. The case
$k=1$ is clear. Recall
that an $A_k$-move $(T,S)$ is a band trivialization of an $A_{k-1}$-move,
say $(T',S')$. Then we have
that $T-S=X_1-X_2$ where $X_1$ and $X_2$ are $1$-singular tangles in
Figure 5.3. Let $Y_1$ and $Y_2$ be $1$-singular tangles in Figure 5.4.
It is clear that $Y_1-Y_2=0$. Thus we have $T-S=(X_1-Y_1)-(X_2-Y_2)$. By
the induction hypothesis we
see that $S'-T'$ is an element of ${\cal R}_{k-1}(S')$.
Therefore both $X_1-Y_1$ and
$X_2-Y_2$ are elements of ${\cal R}_k(T)$. \fin

\begin{center}
\begin{tabular}{cc}
\includegraphics[trim=0mm 0mm 0mm 0mm, width=.4\linewidth]
{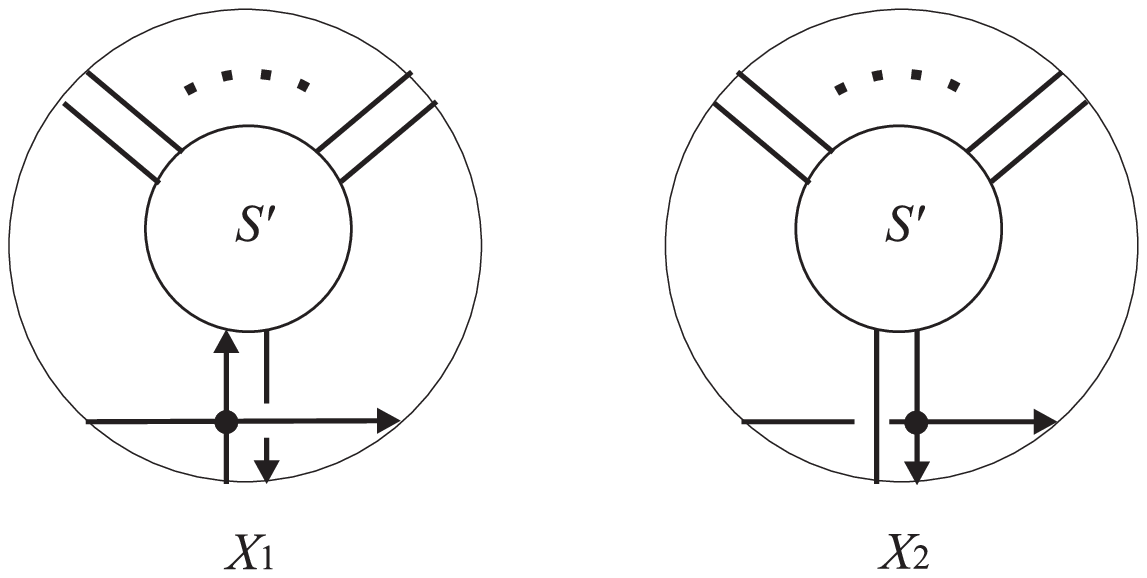}\hspace*{1cm} &
\includegraphics[trim=0mm 0mm 0mm 0mm, width=.4\linewidth]
{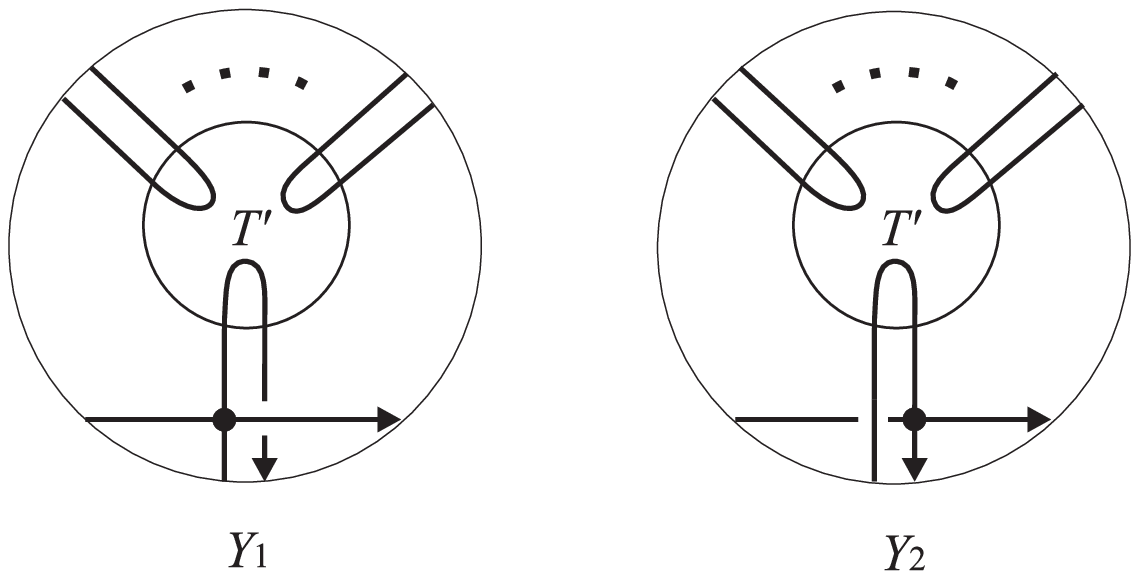}\\
Figure 5.3\hspace*{1cm} & Figure 5.4
\end{tabular}
\end{center}

\medskip\noindent
{\em Proof of Proposition 1.3.} It is sufficient to show that
\[
\sum_{X\subset\{1,2,...,n+1\}}(-1)^{|X|}F(f;\bigcup_{i\in X}\{L_i\},
\bigcup_{i\in X}\{B_i\})
\]
is an element of ${\cal R}_{(n+1)(k-1)}$. We show this together with some
additional claims by
induction on $n$. First consider the case $n=0$. Then we have by Theorem
5.1 and its proof that
$F(f;\emptyset,\emptyset)-F(f;\{L_1\},\{B_1\})$ is a sum of
$(k-1)$-singular embeddings each of
which has all crossing vertices in the associated ball. Note that these
$(k-1)$-singular embeddings
are natural extensions of the $(k-1)$-singular tangles that express the
difference of the
$A_{k-1}$-move, and these $(k-1)$-singular tangles depends only on the link
$L_1$. Next we consider
the general case. Note that
\[\begin{array}{l}
\displaystyle{\sum_{X\subset\{1,2,...,n+1\}}(-1)^{|X|}
F(f;\bigcup_{i\in X}\{L_i\},
\bigcup_{i\in X}\{B_i\})}\\
\hspace*{2em}\displaystyle{=\sum_{X\subset\{1,2,...,n\}}(-1)^{|X|}
F(f;\bigcup_{i\in X}\{L_i\},
\bigcup_{i\in X}\{B_i\})}\\
\hspace*{10em}\displaystyle{-\sum_{X\subset\{1,2,...,n\}}(-1)^{|X|}
F(f;\bigcup_{i\in X}\{L_i\}\cup\{L_{n+1}\},
\bigcup_{i\in X}\{B_i\}\cup\{B_{n+1}\})}.
\end{array}\]
By the hypothesis we have  both
\[
\sum_{X\subset\{1,2,...,n\}}(-1)^{|X|}F(f;\bigcup_{i\in X}\{L_i\},
\bigcup_{i\in X}\{B_i\})
\]
and
\[
\sum_{X\subset\{1,2,...,n\}}(-1)^{|X|}F(f;\bigcup_{i\in
X}\{L_i\}\cup\{L_{n+1}\},
\bigcup_{i\in X}\{B_i\}\cup\{B_{n+1}\})
\]
are sums of $n(k-1)$-singular embeddings and they differ only by the band
sum of $L_{n+1}$.
Therefore we have
\[
\sum_{X\subset\{1,2,...,n+1\}}(-1)^{|X|}F(f;\bigcup_{i\in X}\{L_i\},
\bigcup_{i\in X}\{B_i\})
\]
is a sum of $(n(k-1)+(k-1))$-singular embeddings. \fin

\bigskip
\footnotesize{
 }
\end{document}